\documentclass{elsart3-1}

\usepackage{amssymb}

\usepackage[english,francais]{babel}

\newtheorem{theorem}{Theorem}[section]

\newtheorem{e-proposition}[theorem]{Proposition}

\newtheorem{e-definition}[theorem]{Definition\rm}

\newtheorem{theoreme}{Th\'eor\`eme}[section]

\newtheorem{proposition}[theoreme]{Proposition}
\newtheorem{corollaire}[theoreme]{Corollaire}
\newtheorem{definition}[theoreme]{D\'efinition\rm}

\def\og{\leavevmode\raise.3ex\hbox{$\scriptscriptstyle\langle\!\langle$~}}
\def\fg{\leavevmode\raise.3ex\hbox{~$\!\scriptscriptstyle\,\rangle\!\rangle$}}

\def\K{{\bf K}}\def\N{{\bf N}}\def\P{{\bf P}}\def\Q{{\bf Q}}\def\R{{\bf R}}\def\T{{\bf T}}\def\Z{{\bf Z}}
\def\cA{{\mathcal A}}\def\cH{{\mathcal H}}
\def\arith{{{\mathcal A},\alpha}}\def\Conv{{\rm Conv}}\def\MI{{\rm MI}}\def\init{{\rm init}}\def\irr{{\rm irr}}\def\hnorm{\widehat{h}}

\begin{document}

\begin{frontmatter}




%
\selectlanguage{francais}
\title{G\'eom\'etrie diophantienne et vari\'et\'es toriques}

\vspace{-2.6cm}
\selectlanguage{english}
\title{Diophantine geometry and toric varieties}



\author[authorlabel1]{Patrice Philippon}
\author[authorlabel2]{Mart\'\i n Sombra\thanksref{label2}}

\ead{pph@math.jussieu.fr}
\address[authorlabel1]{Institut de Math\'ematiques de Jussieu, UMR 7586 du CNRS, Case 7012, 175 rue du Chevaleret, 75013 Paris, France.}

\ead{sombra@math.jussieu.fr}
\thanks[label2]{Financ\'e par le programme Ram\'on y Cajal du Minist\`ere Espagnol de la Recherche.}
\address[authorlabel2]{Universitat de Barcelona, Departament d'\`Algebra i Geometria, Gran V\'ia 585, 08007 Barcelona, Espa\~na.}


\begin{abstract}
We present some results on projective toric varieties which are relevant in diophantine geometry. We interpret and study several invariant attached to these varieties in geometrical and combinatorial terms. We also give a B\'ezout theorem for Chow weights of projective varieties an d an application to the theorem of successive algebraic minima. These results are excerpted from~\cite{PSi}\ and~\cite{PSii}. {\it To cite this article: P. Philippon, M. Sombra, C. R. Acad. Sci. Paris, Ser. I  (2005).}

\vskip 0.5\baselineskip

\selectlanguage{francais}
\noindent{\bf R\'esum\'e}
\vskip 0.5\baselineskip
\noindent
Nous pr\'esentons quelques r\'esultats sur les vari\'et\'es toriques projectives, pertinents en g\'eom\'etrie diophantienne. Nous interpr\^etons et \'etudions plusieurs invariants arithm\'etiques attach\'es \`a ces vari\'et\'es en termes g\'eom\'etriques et combinatoires. Nous donnons \'egalement un th\'eor\`eme de B\'ezout pour les poids de Chow des vari\'et\'es projectives et une application au th\'eor\`eme des minimums alg\'ebriques successifs. Ces r\'esultats sont extraits de~\cite{PSi}\ et~\cite{PSii}. {\it Pour citer cet article~: P. Philippon, M. Sombra, C. R. Acad. Sci. Paris, Ser. I  (2005).}

\end{abstract}
\end{frontmatter}

\selectlanguage{english}
\section*{Abridged English version}
For $n,N\in{\bf N}^\times$ let $\cA= ( a_0,  \dots, a_N) \in (\Z^n)^{N+1}$ be a sequence of $N+1$ vectors in $\Z^n$ and consider the diagonal action of the torus ${\bf T}^n:=(\overline{\bf Q}^\times)^n$ on the projective space ${\bf P}^N:={\bf P}^N(\overline{\bf Q})$
$$
*_\cA : \T^n \times \P^N \to \P^N
\quad \quad , \quad \quad  
(s, x) \mapsto (s^{a_0} \, x_0 : \cdots : s^{a_N} \, x_N)
\enspace .
$$ 
The Zariski closure of the orbit of a point $\alpha= (\alpha_0 : \cdots : \alpha_N) \in \P^N$ is denoted $X_\arith := \overline{\T^n *_\cA \alpha} \ \subset \P^N$
and, following~\cite{GKZ}, we call it the {\it projective toric variety} associated to $(\arith)$. This is a subvariety of $\P^N$ stable under the action of $\T^n$, with a dense orbit $X_\arith^\circ:=\T^n *_\cA \alpha$. By restricting the ambiant projective space, we can assume that the point $\alpha$ belongs to the Zariski open subset $({\bf P}^N)^\circ:=\{(x_0:\cdots:x_N):\, x_0\cdots x_N\not=0\}$. In this case the dense orbit $X_{\cA,\alpha^\circ}$ is actually a translate of a subtorus of ${\bf T}^N$, identified to $({\bf P}^N)^\circ$. Without loss of generality we also assume throughout this note that the ${\bf Z}$-module $L_\cA$, generated by the differences of the vectors $a_i$, is equal to $\Z^n$. Under these assumptions the dimension of $X_{\cA,\alpha}$ is equal to $n$.

\smallskip

One can attach several notions of height to a projective variety $X\subset \P^N$. Normalizing with respect to the action of the torus we select the so-called {\it normalized height}, defined by
$$\widehat{h}(X):=\deg(X)\cdot\lim_{k\rightarrow\infty} \frac{1}{k}\frac{h([k]X)}{\deg([k]X)}  \enspace,$$
where $[k]$ denote the $k$-th power morphism $(x_0:\cdots:x_N)\mapsto (x_0^k:\cdots:x_N^k)$ and $h$ any height on the subvarieties of $\P^N$ (the resulting limit is independent of the actual choice of this height, {\it see}~\cite[\S~I.2]{PSi}).

\smallskip

It is well known that, under our hypotheses, the degree of the projective toric variety $X_{\cA,\alpha}$ is equal to $n!$ times the volume of the convex hull $Q_\cA\subset{\bf R}^n$ of the points $a_0,\dots,a_N$, with respect to the usual Lebesgue measure. We obtain an arithmetic analogue of this classical result, expressing the local contributions to the normalized height in a similar way. 

Let $K$ be a numbers field, we denote by $M_K$ the set of absolute values on $K$ extending the $p$-adic and archimedean absolute values of ${\bf Q}$ satisfying $|p|_p=p^{-1}$ and $|2|_\infty=2$. For $v\in M_K$ we denote $K_v$ ({\em resp.} $\Q_v$) the completion of $K$ ({\em resp.} $\Q$) with respect to $v$. Now, for $\alpha\in (K^\times)^{N+1}$ and $v\in M_K$ we consider the polytope $Q_{\cA,\tau_{\alpha\,v}}$, convex hull in $\R^{n+1}$ of the points $(a_0,\log|\alpha_0|_v), \dots, (a_N,\log|\alpha_N|_v)$, and we denote by $\vartheta_{\cA,\tau_{\alpha\,v}}:Q_\cA\rightarrow{\bf R}$ the parametrisation of its upper envelope over $Q_\cA$. Here $\tau_{\alpha\,v}\in{\bf R}^{N+1}$ stands for the vector $(\log|\alpha_0|_v, \dots, \log|\alpha_N|_v)$.  We can now state~:
\smallskip
\begin{theorem}\label{thmi}
{Let $\cA\in({\bf Z}^n)^{N+1}$ such that $L_\cA={\bf Z}^n$ and $\alpha\in (K^\times)^{N+1}$, with the above notations
$$\widehat{h}(X_{\cA,\alpha}) = (n+1)!\sum_{v\in M_K}\frac{[K_v:\Q_v]}{[K:\Q]} \int_{Q_\cA}\vartheta_{\cA,\tau_{\alpha\,v}}(x)dx_1\cdots dx_n \enspace.$$
}
\end{theorem}

\smallskip

The proof requires a sharp normalized arithmetic Hilbert-Samuel theorem and is closely connected to the Chow weights introduced by D.~Mumford in his work~\cite{Mum} on the stability of projective varieties. Furthermore, Theorem~\ref{thmi} has a multiprojective generalization that we present below. We also give a B\'ezout type theorem for the Chow weight of general projective varieties, which in particular implies an exact arithmetic B\'ezout theorem for intersections of toric varieties with monomial divisors.

\selectlanguage{francais}
\section{Poids de Chow des vari\'et\'es toriques}
Soit $\K$ un corps commutatif, $X \subset \P^N(\K)$ une sous-vari{\'e}t{\'e} de dimension $n$ et $\tau= (\tau_0, \dots, \tau_N) \in \R^{N+1}$ un  {\it vecteur  poids}. Consid{\'e}rons la forme de Chow $Ch_X\kern-0.5pt \in \K[U_0, \dots, U_n]$ de $X$ et une variable additionnelle $t$; si l'on {\'e}crit formellement 
$$Ch_X\left( t^{\tau_j} \, U_{i\,j} \, : \ 0\le i \le n, \ 0 \le j \le N \right) = t^{e_0} \, F_0 + \cdots +t^{e_M} \, F_M$$
avec $F_0, \dots, F_M \in \K[U_0, \dots, U_n] \setminus \{ 0\}$ et $e_0 > \cdots > e_M$, le {\it $\tau$-poids de Chow de $X$} est d{\'e}fini par $e_{\tau}(X):= e_0$. Cette notion est introduite dans~\cite[p.~61]{Mum}, (avec des exposants entiers) dans le contexte de la th{\'e}orie g{\'e}om{\'e}trique des invariants et en relation avec l'{\'e}tude de la stabilit{\'e} des vari{\'e}t{\'e}s projectives. 

Pour un vecteur  $\cA= (a_0, \dots, a_N) \in (\Z^n)^{N+1}$ et un poids $\tau \in \R^{N+1}$ on consid{\`e}re le polytope $Q_{\cA,\tau}:= \Conv \left((a_0,\tau_0), \dots, (a_N,\tau_N) \right) \ \subset \R^{n+1}$, dont l'enveloppe sup{\'e}rieure s'envoie bijectivement sur $Q_\cA$ par la projection standard $\R^{n+1} \to \R^n$. Soit alors 
$$\vartheta_{\cA,\tau}  : Q_\cA \to \R \quad , \qquad x \mapsto \max \left\{ y \in \R \, : \ (x,y) \in Q_{\cA,\tau} \right\}\enspace,$$ 
la param{\'e}trisation de cette enveloppe sup{\'e}rieure au-dessus de $Q_\cA$. C'est une fonction  {\it concave} et {\it affine par morceaux}, donc en particulier Riemann int{\'e}grable sur tout polytope. 
\smallskip
\begin{proposition}\label{propi} 
{Soit $\cA \in ( \Z^n)^{N+1}$ tel que $L_\cA= \Z^n$ et $\tau \in \R^{N+1}$, alors 
$$e_{\tau}(X_\cA) = (n+1)! \, \int_{Q_\cA} \, \vartheta_{\cA,\tau}(x) \ dx_1 \cdots  dx_n\enspace.$$}
\end{proposition}

Joint au th\'eor\`eme~\ref{thmi} (dans la version anglaise ci-dessus) cela fournit une interpr\'etation g\'eom\'etrique, en termes de d\'eformations toriques, des contributions locales \`a la hauteur normalis\'ee d'une vari\'et\'e torique.

\section{Th\'eor\`eme de Hilbert-Samuel normalis\'e}
La d{\'e}monstration du th{\'e}or{\`e}me~\ref{thmi} suit une d{\'e}marche indirecte~: au lieu d'utiliser la d{\'e}finition de la hauteur normalis{\'e}e, on s'appuie sur le calcul d'une fonction de Hilbert arith\-m{\'e}\-ti\-que appropri{\'e}e. L'un des principaux obstacles {\`a} surmonter est de trouver une fonction de type Hilbert dont l'expression asymptotique soit li{\'e}e {\`a} la hauteur normalis{\'e}e; les diff{\'e}rentes variantes {\'e}tudi{\'e}es jusqu'{\`a} pr{\'e}sent sont li{\'e}es {\`a} la hauteur projective, {\it voir} par exemple~\cite{GS}, \cite{Ran}. 

Soit $X \subset \P^N$ une vari{\'e}t{\'e} d\'efinie sur le corps de nombres $K$, de dimension $n$, et $I(X) \subset\kern-1pt K[x_0, \dots, x_N]$ son id{\'e}al homog{\`e}ne de d{\'e}finition. Nous utilisons la fonction de Hilbert arithm{\'e}tique ${\cH}_{\rm norm}(X;\cdot)$ associant {\`a} un entier $D$ donn{\'e}, la {\em hauteur de Schmidt} du $K$-espace lin{\'e}aire $I(X)_D$ des formes de degr\'e $D$ de $I(X)$ dans l'espace des formes de degr\'e $D$ de $K[x_0,\dots, x_N]$, identifi\'e \`a $K^{{D+N}\choose N}$ {\it via} la base des mon\^omes \cite[D\'efinition~II.1]{PSi}. Se pose alors la question du comportement asymptotique de cette fonction, {\`a} laquelle nous apportons une r{\'e}ponse pour le cas des vari{\'e}t{\'e}s toriques: 
\smallskip
\begin{proposition}\label{propii} 
{Soit $\cA\in({\bf Z}^n)^{N+1}$ tel que $L_\cA={\bf Z}^n$ et $\alpha\in({\bf P}^N)^\circ$, avec les notations introduites on a 
$$\cH_{\rm norm}(X_\arith;D) = \frac{\widehat{h}(X_\arith)}{(n+1)!} \, D^{n+1} + O(D^n)\enspace.$$ 
lorsque $D$ tend vers l'infini.}
\end{proposition}

Nous faisons remarquer la nature g\'eom\'etrique de cette formule, qui se refl\`ete en particulier dans le comportement en $O(D^n)$ du terme d'approximation, \`a opposer au comportement habituel en $o(D^{n+1})$ dans les th\'eor\`emes de Hilbert-Samuel arithm\'etiques.

\section{Multihauteurs normalis\'ees}
{\`A} l'instar des multidegr\'es, les multihauteurs du tore $\T^n$ plong\'e dans un produit d'espaces projectifs {\it via} plusieurs applications monomiales peuvent aussi s'expliciter \`a l'aide {\it d'int\'egrales mixtes} (ou {\it multi-int\'egrales}) des fonctions concaves apparaissant dans le th\'eor\`eme~\ref{thmi}, {\it voir} la d\'efi\-ni\-tion~\ref{defn} ci-dessous. 
 
Soit $Z\subset\P^{N_0}\times\dots\times\P^{N_m}$ une sous-vari\'et\'e de dimension $n$, on dispose d'une notion de {\it multihauteur projective de $Z$ d'indice $c$} pour chaque $c=(c_0,\dots,c_m)\in\N^{m+1}$ satisfaisant $0\le c_i \le N_i$ et $c_0+\dots+c_m=n+1$. On renvoie \`a~\cite{Rem}\ ou encore~\cite[\S~I.2]{PSi}, pour la d{\'e}finition pr\'ecise. 

Notons  $s:\P^{N_0} \times \cdots \times \P^{N_m} \to \P^{(N_0+1) \cdots (N_m+1) - 1}$ le plongement de Segre, on montre ({\it voir}~\cite[Proposition~I.2]{PSi}) que la suite 
$$  k \mapsto  \deg(s(Z)) \cdot\frac{h_c([k] \, Z)}{k\, \deg([k]\, s(Z))}$$
converge vers une limite $\geq 0$ lorsque $k$ tend vers l'infini, et d\'efinit la {\it multihauteur normalis{\'e}e de $Z$ d'indice $c$}, not\'ee $\hnorm_c(Z)$. 

\medskip

Soient $f: Q \to \R$ et $g: R \to \R$ des fonctions concaves d{\'e}finies sur des ensembles convexes $Q, R \subset \R^n$ respectivement. On pose 
$$f \boxplus g : Q+R \to \R \enspace , \quad \quad x \mapsto \max \{ f(y) + g(z) \, : \ y \in Q, \ z \in R, \ y+z =x\} \enspace,$$
qui est une fonction concave d{\'e}finie sur la somme de Minkowski $Q+R$; on obtient ainsi une structure de semi-groupe commutatif sur l'ensemble des  fonctions concaves (d\'efinies sur des convexes). 
\smallskip
\begin{definition}\label{defn} 
{\rm Pour une famille de $n$ fonctions concaves $f_0 : Q_0 \to \R, \dots,  f_n: Q_n\to \R$ d\'efinies sur des ensembles $Q_i\subset{\bf R}^n$ convexes et compacts, l'{\it int{\'e}grale mixte} (ou {\it multi-int{\'e}grale}) est d{\'e}finie {\it via} la formule
$$\MI(f_0, \dots, f_n) := \sum_{j=0}^n (-1)^{n - j} \kern-3pt \sum_{0 \le i_0 < \cdots < i_j \le n}\kern6pt \int_{Q_{i_0} + \cdots + Q_{i_j}} \kern-3pt f_{i_0} \boxplus \cdots \boxplus 
f_{i_j} \, dx_1 \cdots dx_n \enspace.$$}
\end{definition} 

Cette notion est analogue \`a celle de volume mixte, c'est de plus une fonctionnelle positive,  sym{\'e}trique et lin{\'e}aire en chaque
variable $f_i$, {\it voir}~\cite[\S~IV.3]{PSi}. 

\smallskip

Soit $\cA_0\in(\Z^n)^{N_0+1}$, \dots, $\cA_m\in(\Z^n)^{N_m+1}$ tels que $L_{\cA_0}+\dots+L_{\cA_m}=\Z^n$ et $\underline{\cA}:=(\cA_0,\dots,\cA_m)$. Soit aussi $K$ un corps de nombres, $\alpha_0 \in (K^\times)^{N_0+1}, \dots , \alpha_m \in (K^\times)^{N_m+1}$ et posons $\underline{\alpha}:=(\alpha_0,\dots,\alpha_m)$. Consid\'erons alors l'action monomiale $*_{\underline{\cA}}$ de $\T^n$ sur le produit d'espaces projectifs $\P^{N_0}\times\dots\times\P^{N_m}$ associ\'ee, on note $X_{\underline{\cA},\underline{\alpha}}$ l'adh\'erence de Zariski de l'orbite du point $\underline{\alpha}\in\P^{N_0}\times\dots\times\P^{N_m}$.

Pour chaque $v \in M_K$ on note aussi $\vartheta_{\cA_i,\tau_{\alpha_i\,v}}: Q_{\cA_i} \to \R$  la fonction param{\'e}trant l'enveloppe sup\'erieure du polytope $Q_{\cA_i,\tau_{\alpha_i\,v}} \subset \R^{n+1}$ associ{\'e} au vecteur $\cA_i$ et au poids $\tau_{\alpha_i \, v}=(\log|\alpha_{i\,0}|_v,\dots, \log|\alpha_{i\,N}|_v)$. 
\smallskip
\begin{theoreme}\label{thmii} 
{Soit $c\in\N^{m+1}$ tel que $c_0+\dots+c_m=n+1$, avec les notations ci-dessus on a 
$$\hnorm_{c}(X_{\underline{\cA},\underline{\alpha}}) = \sum_{v\in M_K} \frac{[K_v:\Q_v]}{[K:\Q]} \, \MI_c( \vartheta_{\underline{\cA}, \tau_{\underline{\alpha}\, v}})$$
avec $\MI_c( \vartheta_{\underline{\cA}, \tau_{\underline{\alpha}\, v}}) := \MI\Big( \underbrace{\vartheta_{\cA_0, \tau_{\alpha_0\, v}}, \dots, \vartheta_{\cA_0, \tau_{\alpha_0\, v}}}_{c_0 \hbox{ \rm fois}}
\enspace,\enspace \dots \enspace , \enspace  
\underbrace{\vartheta_{\cA_m, \tau_{\alpha_m\, v}}, \dots, \vartheta_{\cA_m, \tau_{\alpha_m\, v}}}_{c_m \hbox{ \rm fois}} \Big)$.}
\end{theoreme}

\section{Th\'eor\`eme de B\'ezout pour les poids de Chow}
Soit $X\subset{\bf P}^N$ une vari\'et\'e projective quelconque, fixons $\tau \in \Z^{N+1}$, consid{\'e}rons l'action du sous-groupe {\`a} un param{\`e}tre   
$$*_\tau : \T \times \P^N \to \P^N \enspace, \quad \quad (t,(x_0: \cdots : x_N)) \mapsto (t^{\tau_0}\, x_0 : \cdots : t^{\tau_N}\, x_N)$$
et la {\it d\'eformation torique} $ X_\tau$ de $X$ associ{\'e}e, d{\'e}finie comme l'adh\'erence de Zariski de l'ensemble 
$$\{((1:t),  t *_\tau X) \, : \ t\in\T,\ x\in X\}
\subset{\bf P}^1\times{\bf P}^{N} \enspace.$$
La {\it vari{\'e}t{\'e} initiale de $X$ relative au poids  $\tau \in \Z^{N+1}$} est alors d{\'e}finie par
$$\init_\tau(X) := \iota^*(X_\tau \cdot (\{(0:1)\}\times \P^{N}))\enspace,$$
o{\`u} $ \iota: \P^N \to \P^1\times \P^N $ d{\'e}signe l'inclusion $(x_0:\cdots:x_N) \mapsto ((0:1), (x_0:\cdots:x_N))$; c'est donc le {\it cycle limite} $\lim_{t\to \infty} t*_\tau X$ de $X$ sous l'action $*_\tau$, de m{\^e}me dimension et degr{\'e} que $X$.

On montre que lorsque $\tau\in\N^{N+1}$, le poids de Chow de $X$ relatif \`a $\tau$ s'interpr\`ete comme un bi-degr\'e d'une variante de la d\'eformation torique ci-dessus et se comporte donc comme une hauteur. Comme cons\'equence de cette interpr\'etation on obtient dans~\cite[\S~4]{PSii}, un th\'eor\`eme de B\'ezout pour le poids de Chow qui pr\'ecise la majoration obtenue par R.~Ferretti~\cite[Prop.~4.3]{Fer}.
\smallskip
\begin{theoreme}\label{thmiii} {Soit $X\subset\P^N$ une vari\'et\'e projective et $H \in{\rm Div}(\P^N)$ un diviseur ne contenant pas $X$, alors pour $\tau\in\Z^{N+1}$ on a 
\begin{eqnarray*}
e_\tau(X\cdot H) &= &e_\tau(X)\,\deg(H) + e_\tau(H)\,\deg(X) - (\tau_0+\dots+\tau_N)\,\deg(H)\,\deg(X)\cr 
&&\kern5.9cm\displaystyle - \sum_{Y\in\irr(\init_\tau(X))} m(X_\tau\cdot H_\tau;\iota(Y)) \, \deg(Y)
\end{eqnarray*}
o{\`u} la somme porte sur les composantes irr{\'e}ductibles de ${\rm init}_\tau(X)$ et $m(X_\tau\cdot H_\tau;\iota(Y))$ d\'esigne la multiplicit\'e de $\iota(Y)$ dans le cycle intersection $X_\tau\cdot H_\tau$. En particulier, si $H$ est effectif on a pour tout $\tau\in\R^{N+1}$ 
$$e_\tau(X\cdot H) \leq e_\tau(X) \, \deg(H) + e_\tau(H)\,\deg(X) - (\tau_0+\dots+\tau_N)\,\deg(H)\,\deg(X)\enspace,$$
avec \'egalit\'e si et seulement si les vari\'et\'es initiales de $X$ et $H$ s'intersectent proprement.}
\end{theoreme}

\smallskip

On applique ce r\'esultat \`a l'intersection d'une vari\'et\'e torique projective $X_{\cA,\alpha}$ avec un diviseur monomial de la forme ${\rm div}(x^b)$ o\`u $b\in\Z^{N+1}$. Pour chaque place $v\in M_K$ les {\it pans} de la toiture de $Q_{\cA,\tau_{\alpha\,v}}$ induisent, par projection, une {\it d\'ecomposition polyh\'edrale coh\'erente} du polytope $Q_\cA$, {\it voir}~\cite{Stu} par exemple. On note ${\mathcal D}_{\tau_{\alpha\,v}}$ cette d\'ecomposition, dont les \'el\'ements $S$ sont en bijection avec les composantes de la vari\'et\'e initiale $\init_{\tau_{\alpha\,v}}(X_\cA)$. La multiplicit\'e de cette composante dans le cycle intersection $X_{\arith}\cdot {\rm div}(x^b)$ s'\'ecrit alors $[L_{\cA\cap S}:L_\cA]\cdot\vartheta_{\cA\cap S,\tau_{\alpha\,v}}(a)$ avec $a=\sum_{i=0}^Nb_ia_i\in\Z^n$ et $\vartheta_{\cA\cap S,\tau_{\alpha\,v}}$ la fonction param\'etrisant le pan de la toiture de $Q_{\cA,\tau_{\alpha\,v}}$ au-dessus de $S$ ($=Q_{\cA\cap S}$), \'etendue lin\'eairement \`a tout $\R^n$. On obtient donc le th\'eor\`eme de B\'ezout arithm\'etique exact suivant~:
\smallskip
\begin{corollaire}\label{cor} 
Soit $\cA\in (\Z^n)^{N+1}$ tel que $L_\cA=\Z^n$, $\alpha \in (K^\times)^{N+1}$ et $b \in \Z^{N+1}$. Posons  $D:= \sum_{j=1}^N b_j$ et $a:=\sum_{i=0}^Nb_ia_i\in\Z^n$, alors
$$\widehat{h} (X_{\arith}\cdot {\rm div}(x^b)) = D\, \widehat{h}(X_\arith)  - n! \, \sum_{v\in M_K} \frac{[K_v : \Q_v]}{[K:\Q]} \sum_{S\in {\mathcal D}_{\tau_{\alpha v}}} \kern-1pt \vartheta_{\cA\cap S,\tau_{\alpha\,v}}(a){\rm Vol}_{n} (S)\enspace.$$

En particulier, si ${\rm div}(x^b)$ est effectif ({\it c.-\`a-d.} $b\in\N^{N+1}$) on a $\widehat{h} ( X_{\arith}\cdot {\rm div}(x^b) ) \le  D \, \widehat{h}(X_\arith)$.
\end{corollaire}

\section{Optimalit\'e du th{\'e}or{\`e}me des minimums alg\'ebriques successifs} 
\label{minimums}
Soit  $X \subset \P^N$  une vari{\'e}t{\'e} quasi-projective quelconque, 
de dimension $n$.
Le {\em $i$-{\`e}me minimum alg\'ebrique} de $X$ par rapport 
\`a la hauteur normalis\'ee est d{\'e}fini par 
$$
\widehat{\mu}_i(X) := \sup \, \{\ \inf\, \{ \widehat{h}(\xi) \ : \ \xi \in (X \setminus Y)(\overline{\Q}) \} \ : \ Y \subset X, \ {\rm codim}_X (Y ) = i \}
$$
pour $i=1, \dots, n+1$, o{\`u} le supremum est pris sur toutes les
sous-vari{\'e}t{\'e}s $Y$ de codimension $i$ dans $X$.  
On a $\widehat{\mu}_1(X)\geq\dots \geq\widehat{\mu}_{n+1}(X)\geq 0$.

La r{\'e}partition de la  hauteur des points alg{\'e}briques d'une vari\'et\'e {\it ferm\'ee} $X\subset\P^N$ est en relation  avec sa  hauteur, le lien est donn{\'e} par le {\it th{\'e}or{\`e}me des minimums successifs} de
S.-W.~Zhang~\cite[Thm. 5.2 et Lem. 6.5(3)]{Zha95}~:
$$\widehat{\mu}_1 (X)+ \cdots + \widehat{\mu}_{n+1} (X) \, \le \, \frac{\widehat{h}(X) }{\deg (X)} \, \le \,  (n+1) \, \widehat{\mu}_1(X)  
\enspace.$$ 

Comme application du th\'eor\`eme~\ref{thmi}, on construit des exemples toriques montrant qu'\`a des $\varepsilon$-pr\`es, toute configuration possible des minimums successifs se r\'ealise, et que le quotient ${\widehat{h}(X) }/{\deg (X)}$ peut atteindre n'importe quel valeur dans l'intervalle autoris\'e par cet encadrement.
 
\smallskip

\begin{theoreme} \label{densite}
Soient $n, N\in \N$ tels que $N \ge 3\,n+1$ et $\mu_1, \dots, \mu_{n+1}, \nu  \in \R$ tels que $\mu_1\geq \dots\geq \mu_{n+1}\geq 0$ et $\mu_1+\dots +\mu_{n+1} \le \nu < (n+1)\,\mu_1$. Alors pour $0<\varepsilon_1\leq (n+1)\mu_1-\nu$, $\varepsilon_2>0$ arbitraires, il existe une vari\'et\'e torique $X\subset\P^N$ de dimension $n$ telle que 
$$
0 < \mu_i - \widehat{\mu}_i(X) \le \varepsilon_1 \enspace \mbox{ pour } i=1,\dots,n+1 
\mbox{ et } \enspace
\left|\frac{\widehat{h}(X)}{\deg(X)} - \nu\right| < \varepsilon_2\mu_1 \enspace.
$$
De plus, la vari\'et\'e $X$ peut \^etre choisie de degr\'e $\leq (4n^2 \, \varepsilon_2^{-1})^n$ et d\'efinie sur une extension kummerienne $K=\Q(2^{1/\ell})$ de degr\'e $\lfloor\log(2)\,\varepsilon_1^{-1}\rfloor +1$.
\end{theoreme}

\smallskip

\'Etant donn\'e que l'ouvert principal d'une vari\'et\'e torique est le 
translat\'e d'un sous-tore de $(\P^N)^\circ$, il est naturel de comparer cette situation au cas abelien. Soit $A$ une vari\'et\'e abelienne munie d'un fibr\'e ample et sym\'etrique ce qui fournit une notion de hauteur normalis\'ee sur les sous-vari\'et\'es de $A$ (hauteur de N\'eron-Tate pour les points). Soit  $\alpha + B \subset A$ le translat{\'e} d'une sous-vari{\'e}t{\'e} 
abelienne $B$ par un point $\alpha$, et soit ${\rm Tors}(B)$ le sous-groupe des points de torsion de $B$. Si $\beta$ est un point quelconque de $\alpha+B$, alors $\beta+{\rm Tors}(B)$ est un sous-ensemble de points de hauteur $\hnorm(\beta)$ dense dans $\alpha+B$, on en d\'eduit
$$
\widehat{\mu}_1(\alpha +B)= \cdots = \widehat{\mu}_{n+1} (\alpha +B)\enspace. 
$$
Ainsi, dans cette situation l'intervalle du th\'eor\`eme des minimums 
successifs se r\'eduit \`a un point, et on a les \'egalit\'es
$$
\frac{\widehat{h}(\alpha +B)}{\deg(\alpha +B)} 
= (n+1) \, \widehat{\mu}_1 (\alpha +B) 
= \widehat{\mu}_1 (\alpha +B) + \dots + \widehat{\mu}_{n+1} (\alpha +B)\enspace. 
$$

\smallskip

La situation est donc plus riche dans le cas torique, la diff{\'e}rence tenant au fait que les translat{\'e}s de sous-tores de $(\P^N)^\circ$ ne sont pas des ensembles alg\'ebriques ferm{\'e}s. 




\end{document}